\documentclass[12pt, a4]{article}
\usepackage{amssymb}
\def\vol{{\rm vol}}

\def\C{{\bf C}}

\def\bfg{{\bf g}}
\def\bfr{{\bf R}}
\def\SD{{\cal SD}}
\def\ctr{{\rm ctr}}
\def\ghyp{{\bfg_{\rm hyp}}}
\def\gpsi{(\Psi^{-1})^*\bfg}
\def\bfa{{\bf a}}
\def\phia{\phi_{\bf a}}
\def\phiaa{\phi_{(a_0,{\bf a})}}
\def\phihat{\hat{\phi}}
\def\zaa{Z_{(a_0,{\bf a})}}
\def\tzaa{\tilde{Z}_{(a_0,{\bf a})}}
\def\tzbb{\tilde{Z}_{(b_0,{\bf b})}}
\def\xiaa{\xi_{(a_0,{\bf a})}}
\def\div{\mbox{\rm div}}
\def\Mo{\M_1^{\l_0}}

\def\grad{{\rm grad}}
\def\ident{\equiv}
\def\im{{\rm im}}

\def\intersect{\bigcap}
\def\iso{\cong}
\def\lap{\Delta}

\def\lb{\langle}
\def\rb{\rangle}
\def\lnorm{\left\|}
\def\rnorm{\right\|}
\def\na{\nabla}

\def\plus{\oplus}

\def\tensor{\otimes}

\def\vol{{\rm Vol}}

\def\be{\begin{equation}}
\def\ee{\end{equation}}
\def\bearray{\begin{eqnarray}}
\def\eearray{\end{eqnarray}}
\def\bestar{\begin{eqnarray*}}
\def\eestar{\end{eqnarray*}}
\def\ben{\begin{displaymath}}
\def\een{\end{displaymath}}

\def\non{\nonumber}

\def\adp{\mbox {$Ad\: P$}}

\def\A{{\cal A}}

\def\G{{\cal G}}
\def\M{{\cal M}}

\def\lo{\lambda_0}
\def\Mo{{\cal M}_{\lambda_0}}

\newtheorem{theorem}{Theorem}[section]

\newtheorem{lemma}[theorem]{Lemma}

{\addtocounter{theorem}{1}%
{\bf\\[2ex]Remark \arabic{section}.\arabic{theorem} }}{\\[2 ex]}

\def\qed{{\hfill\vrule height10pt width10pt}\\ \vspace{3 ex}}

\def\pf{{\bf Proof}: }

\def\ss{\vspace{.1in}}
\def\bs{\vspace{.2in}}
\def\ms{\vspace{.15in}}
\def\bsn{\vspace{.2in}\noindent}
\def\msn{\vspace{.15in}\noindent}
\def\ssn{\vspace{.1in}\noindent}

\def\a{\alpha}
\def\b{\beta}
\def\d{\delta}
\def\e{\epsilon}

\def\g{\gamma}
\def\i{\iota}

\def\l{\lambda}

\date{12 November 1996}
\title{Instantons and the information metric}
\author{David Groisser\thanks
{Supported in part by National
Science Foundation grant DMS-9307648}
\\ {\small Department of Mathematics}
\\ {\small University of Florida}
\\ {\small Gainesville FL 32611--8105}
\\ {\small USA}
\\ {\small {\it groisser@math.ufl.edu}}
\and Michael K. Murray
\\ {\small Department of Pure Mathematics
}
\\ {\small The University of Adelaide}
\\ {\small Adelaide, SA 5005}
\\ {\small Australia}
\\ {\small {\it mmurray@maths.adelaide.edu.au}}}

\begin{document}

\maketitle

\begin{abstract}
The information metric arises in statistics as a natural inner product
on a space of probability distributions.  In general this inner
product is positive semi-definite but is potentially degenerate.

By associating to an instanton its energy density, we can examine the
information metric {\bf g} on the moduli spaces $\M$ of self-dual
connections over Riemannian 4-manifolds.  Compared with the
more widely known $L^2$ metric, the information metric better reflects
the conformal invariance of the self-dual Yang-Mills equations, and
seems to have better completeness properties. In the case of $SU(2)$
instantons on $S^4$ of charge one, {\bf g} is known to be the
hyperbolic metric on the five-ball.  We show more generally that for
charge-one $SU(2)$ instantons over $1$-connected, positive-definite
manifolds, {\bf g} is nondegenerate and complete in the collar region
of $\M$, and is `asymptotically hyperbolic' there; {\bf g} vanishes at
the cone points of $\M$.  We give explicit formulae for the metric on
the space of instantons of
charge one on $\C P_2$.
\end{abstract}

\setcounter{equation}{0}
\section{Introduction}

The information metric arises in statistics as a metric on a manifold
of probability distributions \cite{Rao}. Its construction is very
simple and can be applied, in principle, to any manifold which
parametrises a set of probability distributions or measures.  However
in this general setting the information metric may be degenerate
(though it is always positive semidefinite).

One class of such manifolds are the minimum sets of variational
problems. Here each point of the minimum set has associated to it an
energy density---a measure of finite total integral---and we can apply
the construction of the information metric on the space of energy
densities. The metric on the minimum set is then
actually the pull-back, under the map that assigns to any point its
energy density, of the information metric on the space of energy
densities. To show that the resulting pull-back metric is
non-degenerate we have to prove that the map sending a point to its
energy density is an immersion. Often in these types of problems there
is a symmetry group which acts and preserves the energy density. The
information metric will clearly be degenerate in directions parallel
to the group action so we have to factor these out and construct a
metric on the minimum set modulo the symmetry group.

For example consider the space of harmonic maps from two-sphere to
itself. Rotations of the target sphere leave the energy density
unchanged so that we consider the quotient space of harmonic maps
modulo rotations.  This is a manifold and it is possible to show that
the energy density is an immersion and that therefore the information
metric is non-degenerate \cite{Murray}.

The space of harmonic maps from $S^2$ to $S^2$ has often been used as
a model for the space of instantons.  In this paper we consider the
problem of showing that the information metric on the instanton moduli
space is non-degenerate and studying its behaviour.  One such example
is already known \cite{Hitchin}.  For instantons of charge $1$ on the
four-sphere the moduli space is a five ball and the information
metric, on general grounds, is conformally invariant and hence the
hyperbolic metric.  Another example is the moduli space of charge one
instantons on $\C P_2$ which is a cone.  We give explicit formulae for
the information metric in this case.  Recall that Donaldson
\cite{Donaldson} showed that for a large class of four manifolds $M$
the moduli space of charge one instantons `interpolates' between these
two models. That is it is a five-dimensional space having  some singular
points which are cones over $\C P_2$, and having an ideal boundary, a
neighbourhood of which looks like $(0, 1) \times M$.  We show that in
such a situation the information metric vanishes at the cone points,
and near the ideal boundary is asymptotic to an ``asymptotically
hyperbolic'' metric of the form ${\rm const}\cdot(dt^2 + g_M) /t^2$
where $g_M$ is the metric on $M$.
We do not know in this generality
if the information metric is nondegenerate between
these two extremes.

The best-known metric on instanton moduli spaces is the $L^2$ metric
([GP1],[GP2],[G1]).  In those cases where nondegeneracy of the
information metric is known, there are several features distinguishing
the the $L^2$ and information metrics on these spaces.  First, unlike
the $L^2$ metric, the information metric is {\em conformally
invariant}, reflecting the conformal invariance of the self-duality
equations (see \S 2).  Second, in the five-dimensional examples above,
and presumably in greater generality, the information metric tends to
be {\em complete} near the ideal boundary (\S 3), whereas the $L^2$
metric tends to be incomplete (indeed the $L^2$ completion is in many
cases known to be the Donaldson/Uhlenbeck compactification; see [F]).
Third, the information metric is truly a quotient metric living on the
unbased moduli space; it is degenerate on the based moduli space (\S
2), and unlike the $L^2$ metric it cannot be induced by a Riemannian
submersion from the based moduli space (\S 5).  Fourth, the
asymptotics of the metrics near cone singularities are quite different
(\S\S 4,5). And fifth, while for the $L^2$ metric one can easily write
down a general formula for the Riemannian connection and curvature in
terms of Green operators, for the information metric there is no
obvious way to write down such general formulas (\S 2).

It was Hitchin \cite{Hitchin} who first suggested the information
metric as an alternative to the $L^2$ metric on these moduli
spaces. The completeness of the information metric in the case of
1--instantons on $S^4$, and its better conformal properties in general
(as compared with the $L^2$ metric), led Hitchin to speculate that for
purpose of differential geometry on instanton moduli spaces the
information metric might be the more suitable of the two.

The remainder of this paper is organized as follows.  In \S 2, we
review the construction of the information metric and see what form it
takes on general instanton moduli spaces.  In \S 3 we restrict
attention to the five-dimensional moduli spaces mentioned above.
We show that the information metric is nondegenerate in the
collar (the region near the ideal boundary) and establish the
asymptotics of the metric there.  In \S 4 we show that for
general $SU(2)$ moduli spaces the metric vanishes at reducible
self-dual connections.  In \S 5 we derive a concrete formula for the
information metric on the space of charge-one instantons on $\C P_2$.
Along the way we discuss the differences noted above between the $L^2$
and information metrics.

Readers interested in the role played by the information metric and
differential geometry in statistics more generally should look at
Amari \cite{Amari} and Murray and Rice \cite{MurrayRice} and
references therein.

\setcounter{equation}{0}
\section{The information metric on instanton moduli space}

   We start by reviewing the definition and some general properties of the
information metric.

   Let $M$ be a compact oriented $n$-dimensional manifold, let
$\Omega^n(M)$ be the space of all (smooth) $n$-forms, and consider the
subspace $\Omega^n_{+}(M)$ of nonnegative $n$-forms that vanish on no
open set.  The open set $\Omega_{++}$ of strictly positive $n$-forms
can be can be viewed as an infinite-dimensional manifold whose
(formal) tangent space at any point $\tau$ is naturally isomorphic to
$\Omega^n(M)$. For any $\tau\in\Omega_{++}$ and any $\a\in
\Omega^n(M)$, the ratio $\a/\tau$ is well-defined, so we obtain
a Riemannian metric on $\Omega_{++}$ by setting
\be\label{grandpa}
{\bfg}_{\rm info}(\a,\b)_\tau=\int_M ({\a}/{\tau}) ({\b}/{\tau})
\tau, \ \ \
\tau\in\Omega_{++}, \ \ \a,\b\in T_\tau\Omega_{++}\iso \Omega^n(M).
\ee
This construction can be generalized by dropping some positivity and
smoothness conditions.  At any $\tau\in\Omega^n_+(M)$ consider the
subspace $H_\tau\subset
\Omega^n(M)$ of all $n$-forms that are of the form $f\tau$ for $f$
some function which is square integrable with respect to $\tau$.
There is an inner product on $H_\tau$ defined by
\be\label{grandma}
\bfg_{\rm info}( f\tau, g\tau )_\tau = \int_M f g \ \tau.
\ee
(More generality is possible but unnecessary for our purposes.)

   We refer to the Riemannian metric (\ref{grandpa}) and its
generalization (\ref{grandma}) as the {\em information metric} on the
space of ``energy densities''.  Note that this metric is invariant
under the group of orientation-preserving diffeomorphisms of $M$.

   Given another manifold ${\cal S}$ and a smooth map
$e:{\cal S}\to\Omega^n_+(M)$ we define the
information metric on ${\cal S}$ by
\ben
\bfg=e^*\bfg_{\rm info}.
\een
The term ``metric'' is used loosely here, since at each $A\in{\cal S}$ the
quadratic form $\bfg_A$ is positive semidefinite but is potentially
degenerate.  For example, if a Lie group $\G$ acts on ${\cal S}$ and leaves
$e$ invariant, then for any vector $\eta\in T{\cal S}$
 tangent to a $\G$-orbit
we have $\bfg(\eta, \cdot)\ident 0$.  When the quotient space
$\M:={\cal S}/\G$ is a manifold it is more appropriate to consider the
induced metric (which we will still denote $\bfg$) on $\M$---although
here too $\bfg$ is potentially degenerate.

    To specialize to gauge theory, let $(M,g_M)$ be a compact,
oriented Riemannian four-manifold and let $P$ be a principal bundle
over $M$ with compact semisimple structure group $G$.  Let $\A$ be the
space of smooth connections on $P$, $\G$ the gauge group,
$\SD\subset\A$ the subspace of self-dual connections, and
$\M=\M(P)=\SD/\G$ be the moduli space of $P$-instantons. ($\M(P)$
depends on $g_M$ but we suppress this from the notation.)  For any
connection $A$ we can define the energy density
\be
e(A) =  ( F_A \wedge \star F_A ) = ( F_A , F_A ) \vol_M
\ee
where $(\ ,\ )$ denotes both the invariant inner product on the Lie
algebra and the inner product constructed out of a metric on $M$. The
volume form $\vol_M$ is defined using the metric $g_M$.

Notice that $e$ is invariant under the action of the gauge group, so
its restriction to $\SD$ induces a smooth map
\be\label{e2}
e \colon \M(P) \to \Omega^4 (M)
\ee
from the moduli space to the space of all four-forms on $M$.

Next recall the definition of the tangent space to $\M(P)$ at
$[A]$. Let $\adp$ be the adjoint bundle. For each self-dual connection
$A$ there is an elliptic complex
\be\label{fes}
0  \to \Omega^0(\adp) {\buildrel d_A \over \longrightarrow}
 \Omega^1(\adp) {\buildrel p_-d_A \over \longrightarrow}
\Omega^2_-(\adp) \to 0
\ee
where $d_A$ is covariant exterior derivative and where $p_-$ projects
a two-form onto its anti-self dual part.  In this complex,
$\Omega^0(\adp)$ represents the Lie algebra of $\G$, $\Omega^1(\adp)$
the tangent space $T_A\A$, and $\Omega^2_-(\adp)$ the tangent space
$T_A\SD$ (at least formally).  The sequence (\ref{fes}) is left-exact
at every irreducible connection, and for generic metrics $g_M$ (and
for some very non-generic metrics) is right-exact at every self-dual
connection.  We will always assume that $g_M$ is such a metric.  In
this case $\M^*(P)$, the subspace of irreducible instantons, is a
manifold, and the tangent space $T_{[A]}\M^*(P)$ is formally $H^1$ of
the complex (\ref{fes}).  We call an object defined on $\M(P)$ smooth
(e.g. the map $e$ in (\ref{e2})) if its restriction to $\M^*(P)$ is
smooth.

The curvature of an instanton cannot vanish on an open set (see the
proof of Theorem 3.4 in \cite{FreedUhlenbeck}), and hence the energy
density of an instanton is contained in the subspace $\Omega^4_+(M)$.
If the image of the tangent space $T_{[A]}\M$ at the energy density to
the energy density $e(A)$ is in the space $H_{e(A)}$ then we can pull
the information metric back to $\SD$.  To see that it is consider the
derivative of the energy density function in the direction of a
tangent vector $\eta$;
\ben{\partial e \over \partial \eta}(A)  =
2 ( d_A\eta \wedge \star  F_A )
=  2( d_A\eta  ,  F_A ) \vol_M
\een
Then
\ben ( d_A\eta  \wedge \star F_A ) =
2{ ( d_A\eta ,  F_A ) \over
( F_A, F_A )} e(A)
\een
and by Cauchy's inequality
\ben
{  ( d_A\eta , F_A ) \over
( F_A, F_A )} \leq
{|d_A\eta| \over |F_A|}
\een
which is certainly square integrable with respect to $e(A)$.

The information metric on $\SD$ (actually on all of $\A$) is therefore
\bearray\non
\bfg(\eta_1, \eta_2) &=& \int_M
\left(\frac{\partial e}{ \partial \eta_1}(A)
/e(A)\right)\left(
\frac{\partial e}{ \partial \eta_2}(A)/e(A)\right) e(A) \\
\label{metric1}
      & =& 4\int_M  {( d_A\eta_1 ,   F_A )
               ( d_A\eta_2  ,  F_A )
                    \over  ( F_A,   F_A ) } \vol_M \\
\non &=& 4\int_M (d_A\eta_1, \hat{F}_A) (d_A\eta_2,\hat{F}_A)
               \vol_M,
\eearray
where $\hat{F}_A=F_A/|F_A|$.
The degeneracy of the information metric tangent to orbits of the
gauge group can be seen explicitly in (\ref{metric1}), since for $v\in
\Omega^0(\adp)$ we have
\be\label{degen}
(F_A,d_A(\eta+d_Av))=(F_A,d_A\eta)+(F_A,[F_A,v])=(F_A,d_A\eta).
\ee
Because of this degeneracy the same formula (\ref{metric1}) serves as
the definition of the quotient information metric on $\M^*(P)$ at
$[A]$.

To have any hope of obtaining a nondegenerate metric on a moduli
space, note that it is important that we divide by the action of the
full gauge group; the information metric on the based moduli space is
automatically degenerate. For if we divide only by the action of the
based gauge group, there is a residual action of $G/({\rm
center}(G))$ on the quotient---which as we saw above is a recipe for
degeneracy along the orbits.

    The information metric is still potentially degenerate on
$\M^*(P)$; conceivably there exist $\eta\notin\im(d_A)$ for which
$(d_A\eta,F_A)\ident 0$.  Since this amounts to an infinite number of
conditions on an element of the finite-dimensional space
$T_{[A]}\M\iso \ker(d_A^-)/\im(d_A)$ (where $d_A^-=p_-d_A$), it is
plausible that for generic metrics $g_M$ the information metric on
$\M^*(P)$ is nondegenerate.  However, any such theorem will need to
make use of irreducibility, since as we show in \S 4 the information
metric vanishes at reducible connections.

Notice that the definition of the information metric does not rely on
choosing a particular base metric $g_M$ except in as much as that is
needed for the definition of self-duality; a conformal change of $g_M$
leaves $\bfg$ unchanged. Moreover it is well known that the space of
instantons is invariant under conformal diffeomorphisms of $M$
homotopic to the identity, and it follows that the information metric
is invariant under transformations of $\M(P)$ induced by conformal
transformations of $M$.  This is in marked contrast to the $L^2$
metric on $\M(P)$ which is only invariant under isometries of
$M$. When $M$ is the four-sphere and $P$ is the $SU(2)$-bundle of
instanton number 1, $\M(P)$ is well-known to be the five-ball, and (as
first noted by Hitchin \cite{Hitchin}) the conformal invariance of the
information metric implies that it is the hyperbolic metric on the
five ball, up to scale.

With a formula for the general information metric in hand, it is
natural to try to write down a formula for the Riemannian connection,
and from this compute curvature.  Indeed, if $\na$ is the putative
Levi-Civita connection, it is straightforward to write down a
formula for $\bfg(\na_\a \b, \g)$, where $
\a,\b,\g$ are vector fields on
$\SD$ obtained by applying the $L^2$ projection $T_A\A\to
\ker((d_A^-)^*)\iso T_A\SD$ to ``constant'' vector fields.  We find
\be\label{lc1}
\left .\bfg(\na_\a \b,\g)\right|_A =
\int_M Q(\a,\b) (d_A\g,\hat{F}_A) \vol_M,
\ee
where
\bestar
Q(\a,\b) &= &
([\a,\b]-d_A (d_A^-)^* (d_A^-(d_A^-)^* )^{-1} p_-[\a,\b], \hat{F}_A) \\
&& +|F_A|^{-1}( (d_A \a,d_A \b) - (d_A \a, \hat{F}_A) (d_A \b,\hat{F}_A) );
\eestar
here $[\a,\b]\in \Omega^2(\adp)$ is the wedge-bracket of
$\a,\b\in\Omega^1(\adp)$.  Were
$Q(\a,\b)$ of the form $(d_A(\mbox {\rm something}),
\hat{F}_A)$ we would obtain a formula for the Levi-Civita
connection---as one is able to do for the $L^2$ metric---but
unfortunately this is not the case here.  For this reason we are at
present unable to analyze curvature invariants of the information
metric, except in those two cases in which the metric is explicitly
computable---the spaces of 1--instantons over $S^4$ and $\C P_2$.

\setcounter{equation}{0}
\section{The collared moduli spaces}

    The best understood moduli spaces are those in which the base
space and principal bundle satisfy the following topological
conditions.

\ssn {\bf Topological Conditions 3.0}
(i) $M$ is a closed, simply connected oriented 4-manifold whose
intersection form on $H^2(M)$ is positive-definite; and
(ii) $P$ is the principal $SU(2)$-bundle over $M$ with instanton
number (Pontryagin index) 1.

\ms    Throughout this section we will assume $M$ and $P$ satisfy these
conditions, which ensure that at smooth points the dimension of
$\M(P)$ is five.  We write $\M(P)=\M_1(M)=\M_1$, and assume $g_M$ has
been chosen so that the subspace $\M_1^*$ of irreducible instantons is
a manifold.  A small neighborhood of each
reducible point in $\M_1$ is then topologically a cone on $\C P_2$; see
[D1] and [FU]. Furthermore there is a compact set in $\M_1$ whose
complement---the ``collar''---consists of instantons with curvature
concentrated near a point, and is diffeomorphic to $(\mbox{\rm open
interval})\times M$ via the map assigning to each sufficiently
concentrated connection $A$ a scale $\l(A)$ and center point
$\ctr(A)$.

    In the explicitly computable examples ($M=S^4$ or $\C P_2$) one
usually defines $\l(A)$ to be the radius of the smallest ball
containing half the total energy $\|F_A\|_2^2$, and $\ctr(A)$ to be
the center of this ball. More generally this definition is modified
non-canonically to ensure differentiability of the scale and
center-point functions (see [D1]); this definition coincides
asymptotically with the preceding one.  In this section, since we deal
with general 4-manifolds obeying the topological conditions 3.0(i), we
use such a non-canonical definition of $\l$ and $\ctr$.  This gives us
a diffeomorphism $\Psi$ from a region $\M_1^{\l_0}$ to $(0,\l_0]\times
M$ for $\l_0$ sufficiently small.

In what follows, we use $\l$ to denote both the scale function on $\M$
and the corresponding real variable in the interval $(0,\l_0]$.  Also
we let $g_M$ denote both the metric on $M$ and its pullback to
$(0,\l_0]\times M$.  The letter $c$ is used for a continually updated
constant whose value can depend on $g_M$ but is independent of all
other parameters of interest.

We will prove the following theorem.

\begin{theorem}\label{ahypthm}
Let $M,P$ satisfy the topological conditions {\rm 3.0}.
Then as $\l\to 0$,
\be\label{ahypthm1a}
\gpsi\sim \frac{128\pi^2}{5}(\frac{d\l^2 +
g_M}{\l^2}):=\ghyp.
\ee
(Here $\sim$ means asymptotic in a $C^0$ sense only: for any $\e>0$,
there exists $\l_0>0$ such that for any tangent vector $X\in
T((0,\l_0]\times M)$ we have $|(\Psi^{-1})^*\bfg (X,X)-\ghyp(X,X)|\leq
\e\ \ghyp (X,X)$.) In particular, for $\l_0$ sufficiently small $\bfg$
is nondegenerate and $\M_1^{\l_0}$ is complete with respect to the
distance function defined by $\bfg$.
\end{theorem}

{\bf Remark.}  The metric $(d\l^2 + g_M)/\l^2$ on $\bfr^+\times M$ is
``asymptotically hyperbolic'' in the sense that as $\l\to 0$, all
sectional curvatures approach 1, regardless of the metric $g_M$.
Also, sitting above every geodesic $\g$ in $M$ is an immersed, totally
geodesic copy of the hyperbolic plane; if $t$ is an arclength
parameter along $\g$, the induced metric on $\bfr^+\times {\bf R}\to
\bfr^+\times {\rm image}(\g)\subset \bfr^+\times M$ is $(d\l^2 +
dt^2)/\l^2$.  Every geodesic in $\bfr^+\times M$ is contained in such
a 2-strip, so $\bfr^+\times M$ is geodesically complete in this
metric.  It follows that for all $a\in \bfr^+$, $(0,a]\times M$ is
complete as a metric space with the induced distance function.  Since
(\ref{ahypthm1a}) implies that the distance functions induced by
$(d\l^2 + g_M)/\l^2$ and $\bfg$ are equivalent on $(0,\l_0]\times M$,
the completeness assertion in Theorem \ref{ahypthm} is an immediate
consequence.

\ss    The proof of Theorem \ref{ahypthm} requires some analysis, to
which we devote the rest of this section.  In several places we use
estimates derived in [GP2], [G2], and [G3].

    Given $[A]\in \M_1^*$, $T_{[A]}\M_1^*$ can be naturally identified
with the harmonic space $H_A= \ker (d_A^*)\intersect \ker
(d_A^-)\subset \Omega^1(\adp)$.  More precisely, the harmonic spaces
piece together into a gauge-invariant subbundle of $\adp$, and
$T_{[A]}\M_1^*$ is naturally isomorphic to the space of
gauge-invariant sections along the gauge orbit $[A]$.  Since the
information metric on $\adp$ is gauge-invariant, to compute the metric
at $[A]$ from (\ref{metric1}) it suffices to choose any $A\in [A]$ and
take $\eta_i$ to be $A$-harmonic.

To make the inner product effectively computable one needs some idea
of what such $\eta_i$ look like.  In [GP2] an ``approximate tangent
space'' $\tilde{H}_A$ was introduced for this purpose, approximating
elements of $H_A$ by purely local objects.  The accuracy of this
approximation in various norms was measured in [GP1] and later
sharpened in [G2-G3], and we will use these estimates to determine the
error introduced by replacing true tangent vectors in (\ref{metric1})
by their $L^2$ projections to $\tilde{H}_A$.  First we recall the
definition of the approximate tangent space.  For this purpose we fix
a smooth cutoff function $b\in C^\infty_0(\bfr)$ with $b(t)=1$ for
$0\leq t\leq 1$, $b(t)=0$ for $t\geq 2$ and $0\leq b(t)\leq 1$
everywhere.  Given $p\in M$, let $r_p$ denote distance to $p$ and
define $\b_p(\cdot)=b(r_p(\cdot)/r_0)$, where $4r_0$ is less than the
injectivity radius of $(M,g_M)$.

\ss
{\bf Definition}.   Let $[A]\in \Mo$ have center point
$p=\ctr(A)$, and let $\{x^i\}_1^4$ be normal coordinates based at $p$.
Given $\bfa\in T_pM$ and $a_0\in \bfr$, define functions
$\phihat=\frac12 \b_p r_p^2, \phia = \b_p a_i x^i$ (note that $\phia$
is independent of the choice of normal coordinate system), and
$\phiaa= \l^{-1}a_0\phihat +\phia$; also define a vector field
$\zaa = \grad(\phiaa)$.  For any vector field $Z$ on $M$ define
$\tilde{Z}^A = \iota_Z F_A$.  The {\em approximate tangent space} at
$A$ is the space

\be\label{def_approxts}
\tilde{H}_A:=\{ \tzaa \mid
(a_0,\bfa)\in \bfr\times T_{p(A)}M \}.
\ee

For $\l_0$ sufficiently small, the $L^2$-orthogonal projection $\pi_A:
\tilde{H}_A\to H_A$  is an isomorphism ([GP2] \S 5), and thus so is
the map
\bearray\non
\a_{\Psi(A)}: \bfr\times T_{\ctr(A)} M &\to & H_A\\
(a_0,\bfa) &\mapsto& -\pi_A\tzaa.
\label{defalpha}
\eearray
In fact, $\a$ approximately inverts the differential of $\Psi$:
\be\label{approxinv}
|(\Psi_*\circ \a_{\Psi(A)} - {\rm Id})(a_0,\bfa)|\leq c
(|a_0|+|\bfa|)\l.
\ee
(The bound $c(|a_0|+|\bfa|)\l^{1-\d}$ appearing in Proposition 5.2 of
[GP2] was strengthened to $c(|a_0|\l^{1+\d}+|\bfa|\l^2)$ in Proposition
1.1 of [G3], but actually any positive power of $\l$ is sufficient for
our purposes.)

    We will prove Theorem \ref{ahypthm} by showing that both $\gpsi$ and
$\ghyp$ are asymptotic to $\a^*\bfg$.

\bsn
\underline{\bf Part I: \ \ \ \ $\a^*\bfg\sim\ghyp$.}

\ss Let $\e>0$ be arbitrary.  We subdivide our argument further into two
steps: computing (\ref{metric1}) when $\eta_i$ are replaced by their
approximate counterparts $\tzaa$, and bounding the error introduced by
the approximation $\tzaa\approx\pi_A\tzaa=-\a_{\Psi(A)}(a_0,\bfa)$.
Below, we write simply $\a$ for $\a_{\Psi(A)}$.

For the first step, it suffices to consider
\be\label{metric2}
\bfg(\tzaa, \tzaa) =
     4 \int_M  \frac{( F_A, d_A\tzaa)^2}{|F_A|^2} \vol_M
\ee
since we can determine $\bfg(\tzaa, \tzbb)$ from this by polarization.

\begin{lemma}\label{lemma1}
For any vector field $Z$ on $M$ and any self-dual connection $A$ we
have
\be\label{lemma1a}
(F_A, d_A(\i_X F_A)) = \frac{1}{2}(\div(X) |F_A|^2 + X(|F_A|^2)).
\ee
Here $\div(X)= -d^*(X^{\rm dual})$, where $X^{\rm dual}$ is the metric
dual of $X$.
\end{lemma}

\pf Let $F=F_A$ let $\{ e_i \}$ be a local orthonormal basis of
$TM$, and let $\theta^i$ be the dual coframe. Write $\na$ for the
Levi-Civita connection on $M$, $\na^A$ for the tensor product
connection on $\adp\tensor \Lambda^k T^*M$ and $\na_i
=\na_{e_i}$. Using the Bianchi identity as in the proof of Lemma 3.1
of [GP2] one finds
$d_A(\i_X F) = \theta^i\wedge (\i_{\na_i X}F) + \na_X^A F.$
Furthermore since $F$ is self-dual,
\ben
(F, \theta^i\wedge (\i_{\na_i X}F)) =
(\i_{e_i}F, \i_{\na_i X}F)
=\frac{1}{2} (e_i, \na_i X)(F,F)
\een
(see [GP2], Lemma 3.4). But $(e_i, \na_i X) =\div(X)$ and $(F, \na_X^A
F)=\frac{1}{2}X(|F|^2)$, so (\ref{lemma1a}) follows.
\qed

   We will apply this with $X=\zaa$ and $A\in
\Mo$.  In view of (\ref{metric2}) we have
\be\label{basic1}
\bfg(\tzaa,\tzaa)= \int_M\left( \div(\zaa)|F_A| + 2\zaa(|F_A|) \right)^2
\vol_M
\ee
We break this up into an integral over $B_{N\l}(p)$ (the ball of
radius $N\l$ centered at $p=\ctr(A)$) and its complement, where $N$ is to
be determined later. (The smaller the $\e$ in Theorem \ref{ahypthm},
the larger $N$ must be taken; to simplify estimates, we will always
take $N\geq 1$.)

\msn\underline{Interior estimates.}

\ss
Recall that for any $k\geq 0$ and any compact
set $K\subset \bfr^4$, for $\l$ sufficiently small, after a suitable
gauge choice and rescaling of coordinates $F_A$ is $C^k(K)$-close to
the curvature of standard instanton of scale 1 on $\bfr^4$.  Ater
undoing the rescalings, this implies that given $N,\d>0$, there exists
$\l_0>0$ such that for all $[A]\in \Mo$ and $x\in B_{N\l}(p)$ we have
\be\label{facta}
\l^2\left| \ |F_A| - |F_{0,\l}| \ \right|(x)
+ \l^3 \left| \ \na |F_A| - \na |F_{0,\l}| \ \right|(x) \leq \d
\ee
(see Theorem 16 of [D1]). Here, after a normal-coordinate
identification of a small ball centered at $p\in M$ with a small ball
centered at $0\in\bfr^4$, $|F_{0,\l}|=\sqrt{48}\l^2/(\l^2+r^2)^2$ is
the norm of the standard instanton on $\bfr^4$ of scale $\l$ centered
at $p$.  It is immaterial in (\ref{facta}) whether the norms are
computed with respect to the metric $g_M$ or the flat metric in normal
coordinates.

For the divergence term in (\ref{basic1}), letting $r=r_p$ (with
$p=\ctr(A)$) we have
\be\label{divz}
\div(\zaa)=-\lap(\phiaa)=4a_0\l^{-1} +O(|\bfa|r+|a_0|\l^{-1}r^2).
\ee
Choose $\l_0$ small enough that $N\l_0\ll r_0$; thus $N\leq {\rm
const}\ c\l^{-1}$ and on $B_{N\l}(p)$ we have $r\leq N\l\leq {\rm
const}$, a fact we will use frequently below without further mention.
Then the derivative of the cutoff $\b$ in the definition of $\zaa$
vanishes on $B_{N\l}(p)$, so, writing $c_1=\sqrt{48}$ and noting that
$|\zaa|\leq c(|\bfa|+|a_0|\l^{-1}r)$ we have
\bestar
\zaa(|F_{0,\l}|)&=& -4c_1\l^2(\l^2+r^2)^{-3}(a_ix^i+a_0\l^{-1}r^2
+O(|\bfa|r^3))\\  &=& -4c_1\l^2(\l^2+r^2)^{-3}(a_ix^i+a_0\l^{-1}r^2)
+O(|\bfa|\l^{-1}).
\eestar
Using (\ref{facta}) and letting $Z=\zaa$, on
$B_{N\l}(p)$ we then have
\bearray\non
\div(Z)|F_A| + 2Z(|F_A|) &=&
\left( 4a_0\l^{-1}|F_{0,\l}| + 2Z(|F_{0,\l}|) \right) \\ \non
&& +O(|\bfa|r+|a_0|\l^{-1}r^2)|F_A|
+4a_0\l^{-1} O(\d\l^{-2}) + O(|Z|\d\l^{-3})\\ \non
&=&
\frac{4c_1\l^2}{(\l^2+r^2)^3}\left( a_0\l^{-1}(\l^2-r^2) -2a_ix^i
\right)\\ &&
+ O\left( (|a_0|+|\bfa|)\l^{-1}(1+N\d\l^{-2})\right) .\non \\
\label{ahypthm2}
\eearray
Because the metric on $B_{N\l}(p)$ is Euclidean (in normal
coordinates) up to $O(r^2)$, we can estimate the main term in
(\ref{metric2}) arising from (\ref{ahypthm2}) by the corresponding
Euclidean integral, and similarly we can bound the integrated error
terms. Writing $d^4x$ for the Euclidean volume form, we find
\bestar
\lefteqn{\int_{B_{N\l}(p)}
\left[ \frac{4c_1\l^2}{(\l^2+r^2)^3}\left( a_0\l^{-1}(\l^2-r^2) -2a_ix^i
\right) \right]^2 d^4x} &&\\
&=& 16c_1^2 {\rm Vol}(S^3)\int_0^{N\l} \left\{
\frac{\l^4}{(\l^2+r^2)^6}\left[
a_0^2\l^{-2}(\l^2-r^2)^2+4\cdot\frac{1}{4} |\bfa|^2r^2+\mbox{\rm (odd
function)} \right] \right\} r^3 dr \\ &=&
16\cdot 48\cdot 2\pi^2 \l^{-2} \int_0^N \frac{\rho^3}{(1+\rho^2)^6}
\left[ a_0^2(1-\rho^2)^2+|\bfa|^2\rho^2\right] d\rho.
\eestar
As $\int_0^\infty \frac{\rho^3}{(1+\rho^2)^6}(1-\rho^2)^2 d\rho$ and
$\int_0^\infty \frac{\rho^5}{(1+\rho^2)^6} d\rho$ converge (both to
1/60), we can choose $N$ large enough that the integrals from $0$
to $N$ above differ from their limiting values by less than
$\e/(16\cdot 48\cdot 2\pi^2)$. Hence
\be\label{ahypthm3}
\left| \int_{B_{N\l}(p)}
\left[ \frac{4c_1\l^2}{(\l^2+r^2)^3}\left( a_0\l^{-1}(\l^2-r^2) -2a_ix^i
\right) \right]^2 d^4x - \|(a_0,\bfa)\|_{\rm hyp}^2
\right| \leq \e \|(a_0,\bfa)\|_{\rm hyp}^2
\ee
where $\|\cdot\|_{\rm hyp}$ is the norm associated with $\ghyp$.

Now we turn to the integrated error terms (still on $B_{N\l}(p)$).
These arise from two sources: the error term in (\ref{ahypthm2}), and
the $O(r^2)$ difference between $d^4x$ and the Riemannian volume
form on $B_{N\l}(p)$.  Noting that
$\frac{\l^2}{(\l^2+r^2)^3}\left|a_0\l^{-1}(\l^2-r^2) -2a_ix^i
\right|\leq c \frac{\l}{(\l^2+r^2)^2} (|a_0|+|\bfa|)$ one finds that
the error introduced by the difference in volume forms is bounded by
the error term arising from (\ref{ahypthm2}), and hence by
\bearray\non
\lefteqn{c(|a_0|+|\bfa|)^2\int_0^{N\l}
\l^{-1}(1+N\d\l^{-2})(\frac{\l}{(\l^2+r^2)^2} +
\l^{-1}(1+N\d\l^{-2}) ) r^3 dr} &&
\\
&\leq & c\l^{-2}(|a_0|^2+|\bfa|^2)(\l^2+N\d)(\log N + \l^2N^4 +\d
N^5).
\label{ahypthm4}
\eearray
Since $\log N\leq |\log\l| + {\rm const}$, by taking $\d=\d(N)$ small
enough (and reducing $\l_0$, if necessary), we can arrange for this last
bound to be less than $\e\|(a_0,\bfa)\|_{\rm hyp}^2$.
Combining this with (\ref{ahypthm3}), we arrive at
\be\label{intestfinal}
\left| \int_{B_{N\l}(p)} \left( \div(\zaa)|F_A| + 2\zaa(|F_A|) \right)^2
\vol_M
- \|(a_0,\bfa)\|_{\rm hyp}^2
\right| \leq 2\e \|(a_0,\bfa)\|_{\rm hyp}^2.
\ee

\msn\underline{Exterior estimates.}

\ss
Let $\Omega=\Omega(p)=B_{2r_0}(p)-B_{N\l(p)}$. Since
\be\label{ext1}
|\div(Z)F_A + 2Z(|F_A|)|
\leq c\left( |\bfa|(r|F_A|+|\na_A F_A|) +|a_0|\l^{-1}(|F_A|+r|\na_A F_A|)
\right)
\ee
and the vector fields $\zaa$ are supported in $B_{2r_0}(p)$, to bound
the contribution to (\ref{basic1}) from the complement of $B_{N\l}(p)$
it suffices to bound $\lnorm r^k F_A \rnorm_{L^2(\Omega)}$ and $\lnorm
r^k\na_A F_A\rnorm_{L^2(\Omega)}$ for $k=0,1$. The non-derivative
norms are easy to evaluate since for $[A]\in\Mo$ ($\lo$ sufficiently
small) we have the pointwise bound
\be\label{ptwsbound}
|F_A|\leq {\rm const} \frac{\l^2}{(\l^2+r^2)^2}
\ee
everywhere on $M$ (see [GP3], \S 5).
 Because the ratio of ${\rm Vol}_M$ to
$d^4x$ is bounded on $B_{2r_0}(p)$, for $-2<k<2$ we therefore have
\be\label{intrkf}
\lnorm r^kF_A \rnorm_{L^2(\Omega)}
\leq c\left[  \int_{N\l}^{r_0}
 \left\{ \frac{\l^4}{(\l^2+r^2)^4} \right\}
r^{2k+3} dr \right]^{1/2}
\leq c(k)\l^k N^{k-2}.
\ee
Bounding the norms that involve $\na_A F_A$ is less direct because one
does not have an analog of (\ref{ptwsbound}) available.  Instead, we
introduce a cutoff function $\g$ that is identically 1 on
$\Omega$:
\ben
\g(\cdot)=\g_p(\cdot)=
b\left(\frac{r_p(\cdot)}{2r_0}\right)\left(
1-b\left(\frac{r_p(\cdot)}{2N\l}\right)\right).
\een
Thus $\lnorm r^k\na_A F_A\rnorm_{L^2(\Omega)}^2 \leq \lnorm \g
r^k\na_A F_A\rnorm_{2}^2$, where $\|\cdot\|_2=\|\cdot\|_{L^2(M)}$.  If
we integrate the $L^2(M)$-inner product by parts and use the
Weitzenb\"{o}ck identity for self-dual 2-forms (which implies that
$|\na_A^*\na_A F_A|\leq c(|F_A|+|F_A|^2)$), we find
\bearray\non
\lnorm \g r^{k} \na_A F_A\rnorm_2 &\leq & c\left(
\lnorm \g k r^{k-1}
F_A\rnorm_2 +\lnorm d{\g}\ r^{k} F_A\rnorm_2
 +\left[ \int_M \g^2r^{2k}|F_A|^3 \vol_M
\right]^{1/2} \right). \non \label{ahypthm5}
\eearray
The first term on the right is $O(\l^{k-1}N^{k-3})$ as in
(\ref{intrkf}), and the third term is similarly seen to be
$O(\l^{k-1}N^{k-4})$ (if $-2<k<4$).  For the middle term, note that
$|d\g|\leq c((N\l)^{-1}\chi_{in}+\chi_{out})$, where $\chi_{in}$
and $\chi_{out}$ are the characteristic functions of the annuli
$B_{N/\l}(p)-B_{N/(2\l)}(p)$ and $B_{4r_0}(p)-B_{2r_0}(p)$
respectively. Integrating as in (\ref{intrkf}) one then finds that
$\lnorm d{\g}\ r^{k} F_A\rnorm_2\leq
c(\l^{k-1}N^{k-3}+\l^2)\leq c\l^{k-1}N^{k-2}$.
\ben
\lnorm r^k\na_A F_A\rnorm_{L^2(\Omega)} \leq \lnorm \g
r^k\na_A F_A\rnorm_2
\leq c(k)\l^{k-1} N^{k-3}
\een
for $-1<k<3$. We conclude that
\bearray\non
\lnorm \div(Z)F_A + 2Z(|F_A|) \rnorm_{L^2(\Omega)}^2
&\leq& c\left(|\bfa|^2(\lnorm rF_A \rnorm_{L^2(\Omega)}^2 + \lnorm \na_A
F_A\rnorm_{L^2(\Omega)}^2)\right. \\ \non
&&  +\left. |a_0|^2\l^{-2}(\lnorm F_A \rnorm_{L^2(\Omega)}^2
+ \lnorm r\na_A F_A\rnorm_{L^2(\Omega)}^2)\right)\\
&\leq& c(|\bfa|^2+|a_0|^2)\l^{-2}N^{-4}.
\label{ahypthm6}
\eearray
Increasing $N$, if necessary (and correspondingly decreasing $\d(N)$
and $\l_0$), we can therefore ensure that the contribution to
(\ref{basic1}) from the complement of $B_{N\l}(p)$ is less than $\e
\|(a_0,\bfa)\|_{\rm hyp}^2$.  Hence
\be\label{ahypthm7}
\left|\ \bfg(\tzaa,\tzaa)-\|(a_0,\bfa)\|_{\rm hyp}^2\ \right|\leq
3\e\|(a_0,\bfa)\|_{\rm hyp}^2.
\ee

   This completes the first step of Part I. For the second step,
define $\xiaa=\pi_A\tzaa-\tzaa=-\a(a_0,\bfa)-\tzaa$.  Then from
(\ref{metric1}) we have
\bestar
|\bfg(\a(a_0,\bfa),\a(a_0,\bfa)) -\bfg(\tzaa,\tzaa)|
\leq 4(2\| d_A\xiaa\|_{2} \| d_A\tzaa \|_{2} +\| d_A\xiaa\|_{2}^2).
\eestar
Pointwise, $d_A\tzaa$ is bounded by the right-hand side of
(\ref{ext1}), and a simpler version of the analysis above shows that
$\| d_A\tzaa \|_{2}\leq c\l^{-1}(|a_0|+|\bfa|).$ {}From Proposition
5.1 of [G2] we have $\| d_A\xiaa\|_{2}\leq c(|a_0|+|\bfa|)$.
Thus $|\bfg(\a_A(a_0,\bfa),\a_A(a_0,\bfa)) -\bfg(\tzaa,\tzaa)| \leq
c\l\| (a_0,\bfa)\|_{\rm hyp}^2\leq \e\| (a_0,\bfa)\|_{\rm hyp}^2$.
Combining this with (\ref{ahypthm6}), we have $\left|\
\a^*\bfg((a_0,\bfa),(a_0,\bfa))-\|(a_0,\bfa)\|_{\rm hyp}^2\
\right|\leq 4\e\|(a_0,\bfa)\|_{\rm hyp}^2$, and thus
\be\label{ahypthm8}
\a^*\bfg\sim\ghyp.
\ee

\bsn\underline{\bf Part II: \ \ \ \ \ $\a^*\bfg \sim\gpsi$.}

\ss Since $\a$ is an isomorphism we can write
$\a(a_0,\bfa)-\Psi^{-1}_*(a_0,\bfa) = \a
(\hat{a}_0,\hat{\bfa})$ for some $(\hat{a}_0,\hat{\bfa})$.  Then $
(\hat{a}_0,\hat{\bfa})=(\Psi_*\circ\a)^{-1} ((\Psi_*\circ
\a - {\rm Id})(a_0,\bfa))$, so from (\ref{approxinv})
\be\label{ahypthm9}
|\hat{a}_0|+|\hat{\bfa}|\leq c\l(|a_0|+|\bfa|) \leq c \l^2\|
(a_0,\bfa)\|_{\rm hyp}
\ee
Letting $\|\cdot\|_{\bfg}$, $\|\cdot \|_{\a^*\bfg}$ , and $\| \cdot
\|_{\gpsi}$ denote the norms
associated with indicated metrics, we have
\bestar
\left| \lnorm (\a_0,\bfa) \rnorm _{\a^*\bfg}^2 -
\lnorm (\a_0,\bfa)\rnorm_{\gpsi}^2 \right|
&=& \left|\lnorm \a (\a_0,\bfa)\rnorm_{\bfg}^2 -
\lnorm (\a_0,\bfa) \rnorm_{\gpsi}^2 \right| \\
&\leq &
\lnorm \a (a_0,\bfa)-\Psi^{-1}_*(a_0,\bfa)\rnorm_{\bfg} \\
&& \cdot \left(
2\lnorm \a (\a_0,\bfa)\rnorm_{\bfg} +
\lnorm
\a (a_0,\bfa)-\Psi^{-1}_*(a_0,\bfa)\rnorm_{\bfg}\right) \\
&=& \lnorm (\hat{a}_0,\hat{\bfa}) \rnorm_{\a^*\bfg}
\left(
2\lnorm (a_0,\bfa) \rnorm_{\a^*\bfg} +
\lnorm (\hat{a}_0,\hat{\bfa}) \rnorm_{\a^*\bfg}\right).
\eestar
But from (\ref{ahypthm8}) and (\ref{ahypthm9})
we have
\ben
\lnorm (\hat{a}_0,\hat{\bfa}) \rnorm_{\a^*\bfg}\leq
2\lnorm (\hat{a}_0,\hat{\bfa}) \rnorm_{\rm hyp} \leq
c\l^{-1}(|\hat{a}_0|+|\hat{\bfa}|)
\leq c\l \lnorm (a_0,\bfa) \rnorm_{\rm hyp} \leq c
\l \lnorm (a_0,\bfa) \rnorm_{\rm \a^*\bfg},
\een
so $\a^*\bfg\sim \gpsi$, as desired.

This completes the proof of Theorem \ref{ahypthm}. \qed

\setcounter{equation}{0}
\section{The information metric at reducible connections}

The next theorem shows that at reducible $SU(2)$ connections, the
information metric is not merely degenerate, but actually zero.

\begin{theorem}
Let $P$ be a principal $SU(2)$-bundle and assume the base metric $g_M$
is one for which $\SD$ is a manifold.  Let $\bfg$ denote the
information metric on $\SD$. Then at every reducible
self-dual connection, $\bfg=0$.
\end{theorem}

\pf Let $A\in{\cal SD}$ be reducible and let $\eta\in T_A\SD=
\Omega^1(\adp)$; thus
$d_A\eta=\star d_A\eta$. Since $A$ is a reducible $SU(2)$ connection,
there exists a nonzero covariantly constant section
$\Phi\in\Gamma(\adp)$, and moreover $F_A=\Phi\tensor\omega$ for some
real-valued self-dual 2-form $\omega$ (see [FU], Theorem
3.1). Momentarily let $(\cdot,\cdot)$ denote the pairing
$\Omega^k(\adp)\tensor \adp\to \Omega^k(M)$ given by taking inner
product only on the $\adp$ factors, and let
$\a=(\eta,\Phi)\in\Omega^1(M)$.  Since $\Phi$ is covariantly constant,
we obtain $d\a=(d_A\eta, \Phi)=(\star d_A\eta, \Phi)=*d\a$.  Therefore
$d^*d\a= -\star dd\a=0$, implying $\lb d\a,d\a \rb_{L^2} = \lb \a,
d^*d\a\rb_{L^2} =0$, and hence $d\a=0=(d_A\eta,\Phi)$.  Thus
$d_A\eta$ is pointwise perpendicular to $\Phi$, hence to $F_A$, so
from (\ref{metric1}) we have $\bfg(\eta,\cdot)=0$. \qed

{\bf Remarks}. (1) Even without the assumption that $\SD$ is a
manifold, our argument shows that the information metric vanishes on
the formal tangent space to $\SD$.  In contrast, the $L^2$ metric is
positive definite on every subspace of $\Omega^1(\adp)$.

(2) Because of the quadratic nature of the integrand in the definition of
$\bfg$, the analysis above shows that $\bfg$ vanishes to
order two at a reducible connection.

\setcounter{equation}{0}
\section{The information metric on $\M_1(\C P_2)$}

Let $g_{FS}$ be the Fubini-Study metric on $M=\C P_2$, with sectional
curvatures between 1 and 4.  The moduli space $\M_1=\M_1(\C P_2)$ of
$SU(2)$ 1-instantons over $M$ is a cone on $\C P_2$. Let $\l$ be the
usual ``scale function'' on these instantons: $\l(A)$ is the radius of
the smallest ball containing half the total energy $\| F_A\|_2^2$.
The vertex $[A_0]$ of the cone is a reducible (and homogeneous)
connection, and $\l(A_0)=1$, while $\l(A)\to 0$ as $[A]$ approaches
the ``ideal boundary'' of $\M_1$.  For $[A]\in \M_1$ other than
$[A_0]$, the center of the ball defining $\l(A)$ is unique, and the
map sending $[A]$ to its scale and center point is a diffeomorphism
from the punctured cone $\M_1^*=\M_1-\{[A_0]\}$ to $(0,1)\times M$. We
will implicitly use this identification of $\M_1^*\iso (0,1)\times M$
below.

\begin{theorem}\label{cp2thm1}
On the punctured moduli space $\M_1(\C P_2)-\{[A_0]\}$, the information
metric $\bfg$ is given by
\be\label{cp2thm1a}
\bfg=\frac{128}{5}\pi^2\left( \frac{f(\l)d\l^2 + h(\l)g_{FS}}{\l^2} \right),
\ee
where
\bearray
f(\l) &=&\frac{1-\frac{7}{3}\l^2 +\frac{14}{9}\l^4 -\frac{2}{3}\l^6
+\frac{2}{27}\l^8}{1-\l^2}
-\frac{\frac{30}{81}\l^8-\frac{20}
{81}\l^{10}}{(1-\l^2)^2}\log\frac{\l^2}{3-2\l^2}
,\\ \nonumber && \\
h(\l)&=&\frac{1-\frac{7}{3}\l^2 +\frac{23}{18}\l^4 +\frac{93}{108}\l^6
-\frac{77}{108}\l^8}{1-\l^2}
+\frac{\frac{5}{18}\l^6
-\frac{10}{27}\l^8+\frac{1}{81}\l^{10}}{(1-\l^2)^2}
\log\frac{\l^2}{3-2\l^2} .\nonumber \\
\eearray

\end{theorem}

   Before giving the proof, we review the parametrization of $\M_1(\C
P_2)$ discussed in [G1], where the $L^2$ metric $\bfg_2$ was computed.
The $\C^2$-bundle over $\C P^2$ of instanton number 1 is $L\plus
L^{-1}$, where $L$ is the hyperplane bundle.  These are homogeneous
bundles under the action of of $SU(3)$ on $\C P^2$ induced by the
standard linear action on $\C^3$, and $SU(3)$ preserves the space of
self-dual connections. The canonical connection $A_0$ on the
holomorphic hermitian vector bundle $L\plus L^{-1}$ is a fixed point
of this action.  It was shown in [G1] that an identification of $\C^3$
with the formal tangent space $T_{[A_0]}\M_1$ induces an
identification of $\M_1/SU(3)$ with $\C^3/SU(3)\iso [0,\infty)$, and
of any given orbit in $\M_1^*:=\M_1-\{[A_0]\}$ with
$SU(3)/S(U(1)\times U(2))\iso \C P_2$.  It is more convenient to
replace $[0,\infty)$ with $[0,1)$. Also given in [G1] is a 1-parameter
family of self-dual connections $\{A_t\}, 0\leq
t=\sqrt{1-\l^2(A_t)}<1$, centered at $p_0=[1,0,0]\in \C P_2$. The
image of this family in $\M_1$ is transverse to the $SU(3)$-orbits
(i.e. is a section of the fibration $\M_1\to \M_1/SU(3)$).  With these
identifications we can speak of ``tangential'' and ``radial''
directions in $T\M_1^*$, i.e. those tangent to the $SU(3)$-action and
those tangent to this 1-parameter family (or to its translates under
the $SU(3)$-action).  For $t>0$, we will write $\C P_2^t$ for the
$SU(3)$-orbit through $[A_t]$.

\bs
{\bf Proof of Theorem \ref{cp2thm1}}. The explicit formulas in \S 4 of
[G1] for elements of $T_{[A_t]}\M_1$, which we do not repeat here,
give all the data needed to compute inner products; we only state the
relevant consequences below.  The computation for the information
metric is actually simpler than that for the $L^2$ metric $\bfg_2$,
since by (\ref{degen}) there is no need to compute the
$L^2$-orthogonal projection of elements of $T_{A_t}\SD$ to
$\ker(d_{A_t}^*)$ (the difference between $\eta\in\Omega^1(\adp)$ and
its projection is in the image of $d_{A_t}$).  The action of $SU(3)$
on $\M_1$ is preserves $\bfg$ (as well as $\bfg_2$), so it suffices to
determine the metric at each $[A_t]$; the general form
(\ref{cp2thm1a}) then follows from symmetry.

If $\eta_{rad}$ is the ``radial'' tangent vector $dA_t/dt$ (where
$\{A_t\}$ is the family above), one finds from [G1] that pointwise
\ben
(F_{A_t}, d_{A_t}\eta_{rad}) =
\frac{32t(1-t^2)D^3}{(D-t^2)^5}\left( -D^2+D(3-4t^2) + 3t^2-t^4
\right),
\een
where $D=1+|z_1|^2 +|z_2|^2$, and where $(z_1,z_2)$ are the usual
coordinates on a standard $\C^2\subset \C P_2$.  The action of $SU(3)$
induces a complex-linear isometric identification of the dual space
$(\C^2)^*$ (equipped with its standard hermitian metric) with
$T_{p_0}\C P_2$ (equipped with the metric $g_{FS}$), and therefore with
$T_{[A_t]}\C P_2^t$.  If $\eta_{\mu,t}\in T_{[A_t]}\C P_2^t$ denotes
the ``tangential'' tangent vector corresponding to $\mu\in (\C^2)^*$,
then from [G1] one finds that
\ben
(F_{A_t}, d_{A_t}\eta_{\mu,t}) =
\frac{-96t^2(1-t^2)^2D^3(D+t^2)\mbox{\rm Re}(\mu(z_1,z_2))}
{(D-t^2)^5}.
\een

Finally, one has $|F_{A_t}|^2= 16D^3(1-t^2)^2(D+2t^2)(D-t^2)^{-4}$,
and the Riemannian volume form on $\C P^2$ is $D^{-3}d^4x$, where
$d^4x$ is the standard volume form on ${\bf R}^4\iso \C^2$.  The proof
of Theorem \ref{cp2thm1} is now a matter of computing integrals and
substituting $t=\sqrt{1-\l^2}$. \qed

The $L^2$ metric $\bfg_2$ also has the form $f(\l)d\l^2+h(\l)g_{FS}$
(with different coefficient functions $f,h$), and it is interesting to
compare the behavior of the two metrics near the vertex $[A_0]$ of the
cone.  For simplicity, we rescale the metrics as indicated in the
table below.  Near $[A_0]$ it is more natural to express $\bfg$ in
terms of $r$, distance to the vertex, rather than in terms of $\l$.
For each metric, let $N$ denote a unit vector in $T_A\M_1$ normal to
$\C P_2^{t(r)}$, and let $T,T'$ denote unit vectors tangent to such a
$\C P_2$, with $T'$ orthogonal both to $T$ and to $JT$, where $J$ is
the complex structure on $\C P_2$.  All sectional curvatures of $\M_1$
(with respect to either metric) can be expressed in terms of the three
``primary'' sectional curvatures $\sigma_{TN}:=\sigma(T,N)$ (the
curvature of the two-plane spanned by $T$ and $N$),
$\sigma_{TT4}:=\sigma(T,JT)$, and $\sigma_{TT1}:=\sigma(T,T')$.  (In
the Fubini-Study metric $g_{FS}$ on $\C P_2$, these last two sectional
curvatures are 4 and 1, respectively.)  As $r\to 0$, one
has the following asymptotics.

%%%%%%%%%%%%%%%%%%%%%%%%%%%%%%
\begin{table}[h] \label{table1}
%\begin{center}{TITLE?}\end{center}
\begin{tabular}{c|c|c|c}
 metric & $\sigma_{TN}$ & $\sigma_{TT1}$ & $\sigma_{TT4}$ \\
\hline &&&
\\
$(\frac{128\pi^2}{5})^{-1}\bfg_{\rm{info}} =
 dr^2+r^2(3+O(r^2))g_{FS}$ & $-\frac{8}{125} +O(r)$ &
$-\frac{2}{3}r^{-2}+O(1)$&
$\frac{1}{3}r^{-2}+O(1)$ \\ &&& \\
\hline &&& \\
$(4\pi^2)^{-1}\bfg_2 =
 dr^2+r^2(1+O(r^2))g_{FS}$ &
  $-\frac{3}{2}+O(r^2)$ & $-\frac{3}{2}+O(r^2)$ &
$3r^{-2}+O(1)$
\\ &&& \\
\hline
\end{tabular}
\caption{ Comparison of the
information metric $\bfg=\bfg_{\rm info}$ and
the $L^2$ metric $\bfg_2$ on $\M_1(\C P_2)$.}
\end{table}
%%%%%%%%%%%%%%%%%%%%%%%%%%%%%

    There are several interesting observations that can be drawn from
this table.  Recall that $\M_1$ is the quotient by $SO(3)$ of a smooth
8-dimensional manifold $\tilde{\M}_1$.  $SO(3)$ acts freely on the
irreducible connections in $\tilde{\M}_1$, and the stabilizer of each
reducible connection is a circle.  In [GP2] it was shown that there is
a metric on $\tilde{\M}$ for which the quotient map, restricted to
irreducible connections, is a Riemannian submersion.  It was also
shown that given any principal Riemannian submersion of this general
type---i.e. with $SO(3)$ acting smoothly on a Riemannian manifold,
freely except at isolated points stabilized by circles, so that the
singularities in the quotient are cones on $\C P_2$
topologically---the asymptotics of the base metric near the singular
points is {\em always} of the form $dr^2+r^2g_{FS}+O(r^4)$.  The
factor of 3 in the asymptotics of the information metric therefore
shows that {\em there is no smooth metric on $\tilde{\M}_1$ for which
the map $\tilde{\M}_1^*\to \M_1^*$ is a Riemannian submersion}.  It
would be nice to have an interpretation of this ``3'' in in terms of a
model geometry, or family of model geometries, that universally gives
the behavior of the information metric near cone singularities.

While the authors do not know at this writing whether there is such a
universal model, a candidate is the following. On $\bfr^6$, let $\rho$
denote distance to the origin in the standard flat metric, let
$g_{S^5}$ be the standard metric on $S^5$, and using polar coordinates
define $g_{sing}=\rho^2(d\rho^2 + \frac{3}{4}\rho^2g_{S^5})$.  This is a
smooth field of quadratic forms on $\bfr^6$, but vanishes to order two
at the origin, just as the information metric on $\SD$ vanishes to
order two at a reducible connection.  We can still use $g_{sing}$ to
define distance to the origin; if this distance is $r$, then
$r=\rho^2/2$.  On the complement of the origin, $g_{sing}=dr^2+
3g_{S^5}$.  If we now make the identification $\bfr^6\iso\C^3$ and let
$U(1)$ act the usual way, the quotient is a cone on $\C P_2$.  The
circle action preserves $g_{sing}$, and on the complement of the
vertex we obtain a Riemannian submersion metric $dr^2+3g_{FS}$.

    While this model works well for $\M_1(\C P_2)$, an examination of
the formula for $\bfg$ near more general cone singularities gives no
reason to believe that such a symmetric model is valid more generally.

    Finally, we mention that the table also shows qualitative
differences in the sectional curvatures $\bfg_{\rm info}$ and
$\bfg_2$: for $\bfg$, the sectional curvatures are unbounded both
positively and negatively, while for $\bfg_2$ they are only unbounded
positively.  Note also that the difference between $\sigma_{TN}$ and
its limiting value is $O(r^2)$ for $\bfg_2$, but $O(r)$ for $\bfg$.
Again, the generality of these features is unclear.

\bs
\noindent{\bf Acknowledgements:} The financial support of the
National Science Foundation (USA), the Australian Research Council
and the University of Adelaide
are gratefully acknowledged. MKM would like to thank Nick Buchdahl and
Alan Carey for many useful conversations. DG would like to thank the
Pure Maths department of University of Adelaide for its hospitality
during the summer of 1993, when much of this work was completed.

\end{document}